\documentclass{gtart_h}

\def\ifplaintex{\expandafter\ifx\csname documentclass\endcsname\relax}

\def\gtp{{\mathsurround=0pt\it $\cal G\mskip-2mu$eometry \&\ 
$\cal T\!\!$opology $\cal P\!$ublications}}  

\def\recd{{\small Received:\qua\receiveddate\ifx\reviseddate\relax
\else\qquad Revised:\qua\reviseddate\fi\par}} 

\def\lognumber#1{\def\thelognumber{#1}}
\def\volumenumber#1{\def\thevolumenumber{#1}}
\def\volumeyear#1{\def\thevolumeyear{#1}}
\def\papernumber#1{\def\thepapernumber{#1}}
\def\pagenumbers#1#2{\def\startpage{#1}\def\finishpage{#2}}
\def\published#1{\def\publishdate{#1}}

\def\received#1{\def\receiveddate{#1}}
\def\revised#1{\def\reviseddate{#1}}
\def\accepted#1{\def\accepteddate{#1}}

\def\asciiaddress#1{\def\theasciiaddress{#1}}
\def\asciiemail#1{\def\theasciiemail{#1}}

\long\def\asciiabstract#1{\long\def\theasciiabstract{#1}}

\let\\\par\let\thelognumber\relax\let\thevolumenumber\relax
\let\thepapernumber\relax\let\thevolumeyear\relax\let\startpage\relax
\let\finishpage\relax\let\publishdate\relax\let\receiveddate\relax
\let\reviseddate\relax\let\accepteddate\relax\let\theasciititle\relax
\let\theasciiauthors\relax\let\theasciiaddress\relax
\let\theasciiabstract\relax

\let\theasciiemail\relax

\ifplaintex
\font\logobig=cmssbx10 scaled 3836
\font\logomed=cmssbx10 scaled 2557
\else
\font\logobig=cmssbx10 scaled 4200
\font\logomed=cmssbx10 scaled 2800
\fi

\long\def\makeagttitle{   
\count0=\startpage
\agt\hfill      
\hbox to 45truept{\vbox to 0pt{\vglue -13truept{\logomed A\kern -.37em{\logobig 
T}\kern -.38em G}\vss}\hss}
\break
{\small Volume \thevolumenumber\ (\thevolumeyear)
\startpage--\finishpage\nl
Published: \publishdate}

\vglue .25truein

{\parskip=0pt\leftskip 0pt plus
1fil\def\\{\par\smallskip}{\Large\bf\thetitle}\par\medskip} \vglue
0.05truein

{\parskip=0pt\leftskip 0pt plus 1fil\def\\{\par}{\sc\theauthors}
\par\medskip}%
 
\vglue 0.03truein 

{\small\leftskip 25truept\rightskip 25truept{\bf Abstract}\stdspace\theabstract

{\bf AMS Classification}\stdspace\theprimaryclass
\ifx\thesecondaryclass\relax\else; \thesecondaryclass\fi\par
{\bf Keywords}\stdspace \thekeywords\par}\vglue 7truept

}   

\ifplaintex
\font\phead=cmsl9 scaled 950
\font\pnum=cmbx10 scaled 913
\font\pfoot=cmsl9 scaled 950
\headline{\vbox to 0pt{\vskip -4.5mm\line{\small\phead\ifnum
\count0=\startpage ISSN 1472-2739 (on-line) 1472-2747 (printed)
\hfill {\pnum\folio}\else\ifodd\count0\def\\{ }%
\ifx\theshorttitle\relax\thetitle\else\theshorttitle\fi\hfill{\pnum\folio}
\else\def\\{ and }{\pnum\folio}\hfill\ifx\theshortauthors\relax\theauthors
\else\theshortauthors\fi\fi\fi}\vss}}
\footline{\vbox to 0pt{\vglue 0mm\line{\small\pfoot\ifnum\count0=\startpage
\copyright\ \gtp\hfill\else
\agt, Volume \thevolumenumber\ (\thevolumeyear)\hfill\fi}\vss}}
\else
\font=cmsl9 scaled 1050
\font\lnum=cmbx10 
\font=cmsl9 scaled 1050
\makeatletter
\def\@oddhead{{\small\lhead\ifnum\count0=\startpage ISSN 1472-2739 
(on-line) 1472-2747 (printed)\hfill {\lnum\number\count0}\else\ifodd\count0
\def\\{ }\ifx\theshorttitle\relax \thetitle \else\theshorttitle\fi\hfill
{\lnum\number\count0}\else\def\\{ and }{\lnum\number\count0}
\hfill\ifx\theshortauthors\relax 
\theauthors\else\theshortauthors\fi\fi\fi}}\def\@evenhead{\@oddhead}
\def\@oddfoot{\small\lfoot\ifnum\count0=\startpage\copyright\ \gtp\hfill\else
\agt, Volume \thevolumenumber\ (\thevolumeyear)\hfill\fi}
\def\@evenfoot{\@oddfoot}
\makeatother
\fi
\let\maketitlepage\makeagttitle
\let\makeshorttitle\maketitlepage
\let\maketitle\maketitlepage

\newwrite\gtoutfile
\long\gdef\makeheadfile{  
{\def\\{, }\def\s{ }
\immediate\openout\gtoutfile head.xxx
\immediate\write\gtoutfile{Proxy-for: \ifx\theasciiauthors\relax
\theauthors\else\theasciiauthors\fi\s<\ifx\theasciiemail\relax\theemail\else\theasciiemail\fi>}
\immediate\write\gtoutfile{\noexpand\\}
\immediate\write\gtoutfile{Authors: \ifx\theasciiauthors\relax
\theauthors\else\theasciiauthors\fi}
{\def\\{ }\immediate\write\gtoutfile{Title: \ifx\theasciititle\relax
\thetitle\else\theasciititle\fi}}
\immediate\write\gtoutfile{Subj-class: GT or SG, GR etc}
\immediate\write\gtoutfile{MSC-class: \theprimaryclass\ifx\thesecondaryclass\relax\else, \thesecondaryclass\fi}
\immediate\write\gtoutfile{Journal-ref: Algebr. Geom. Topol. \thevolumenumber\s
(\thevolumeyear) \startpage-\finishpage}
\immediate\write\gtoutfile{Comments: Published by Algebraic and
Geometric Topology at}
\immediate\write\gtoutfile{\s\s\s  http://www.maths.warwick.ac.uk/agt/AGTVol\thevolumenumber/agt-\thevolumenumber-\thepapernumber.abs.html}
\immediate\write\gtoutfile{\noexpand\\}
\immediate\write\gtoutfile{}
\ifx\theasciiabstract\relax
\immediate\write\gtoutfile{\theabstract}\else
\immediate\write\gtoutfile{\theasciiabstract}\fi
\immediate\write\gtoutfile{}
\immediate\write\gtoutfile{\noexpand\\}
\immediate\write\gtoutfile{}
\immediate\closeout\gtoutfile}}  

\def\maketitlepage{\makeagttitle\makeheadfile}
\let\makeshorttitle\maketitlepage
\let\maketitle\maketitlepage

\lognumber{38}
\volumenumber{5}
\volumeyear{2005}
\papernumber{38}
\pagenumbers{911}{922}
\received{17 September 2004} 
\revised{30 May 2005}
\accepted{1 June 2005}
\published{3 August 2005}

\usepackage{amssymb,amsmath}

\def	\Z	{{\mathbb Z}}
\def	\R	{{\mathbb R}}
\def	\C	{{\mathbb C}}
\def	\P	{{\mathbb P}}
\def	\CP	{{\C\P}}

\def	\t	{{\mathfrak t}}

\def	\ft	{{\mathfrak t}}

\def	\inv	{^{-1}}

\def    \calT   {{\mathcal T}}

\def	\d	{\operatorname{d}}
\def	\std	{{\operatorname{std}}}

\def	\U	{{\operatorname{U}}}

\def	\codim	{\operatorname{codim}}

\def	\image	{\operatorname{image}}

\def	\Gr	{\operatorname{Gr}}

\def	\reg	{{\operatorname{reg}}}
\def	\to	{\longrightarrow}
\def	\ssminus	{\smallsetminus}

\def	\half	{\frac{1}{2}}

\def	\calJ		{{\mathcal J}}

\numberwithin{equation}{section}
\newtheorem {Theorem}			{Theorem}

\newtheorem {Proposition} [equation]	{Proposition}

\theoremstyle{definition}
\newtheorem{Definition}[equation]{Definition}

\theoremstyle{remark}
\newtheorem{Remark}[equation]{Remark}

\newtheorem{Example}[equation]{Example}

\begin{document}
\title{The Gromov width of complex Grassmannians}

\author{Yael Karshon\\Susan Tolman}
\address{Department of Mathematics, the University of Toronto\\Toronto, 
Ontario, M5S 3G3, Canada\\{\rm and}\\Department of Mathematics, 
University of Illinois at Urbana-Champaign\\1409 W Green St, Urbana, IL 61801, USA}
\asciiaddress{Department of Mathematics, the University of Toronto\\Toronto, 
Ontario, M5S 3G3, Canada\\and\\Department of Mathematics, 
University of Illinois at Urbana-Champaign\\1409 W Green St, Urbana, IL 61801, USA}
\gtemail{\mailto{karshon@math.toronto.edu}{\qua\rm 
and\qua}\mailto{stolman@math.uiuc.edu}}
\asciiemail{karshon@math.toronto.edu\\stolman@math.uiuc.edu}

\begin{abstract}
We show that the Gromov width of the Grassmannian of complex $k$-planes 
in $\C^n$ is equal to one when the symplectic form is normalized 
so that it generates the integral cohomology in degree 2.
We deduce the lower bound from  more general results.
For example, if a compact manifold $N$ with an integral symplectic form  
$\omega$ admits a Hamiltonian circle action with a fixed point $p$
such that all the isotropy weights at $p$ are equal to one,
then the Gromov width of $(N,\omega)$ is at least one.
We use holomorphic techniques to prove the upper bound.
\end{abstract}

\asciiabstract{%
We show that the Gromov width of the Grassmannian of complex k-planes
in C^n is equal to one when the symplectic form is normalized so that
it generates the integral cohomology in degree 2.  We deduce the lower
bound from more general results.  For example, if a compact manifold N
with an integral symplectic form omega admits a Hamiltonian circle
action with a fixed point p such that all the isotropy weights at p
are equal to one, then the Gromov width of (N,omega) is at least one.
We use holomorphic techniques to prove the upper bound.}

\primaryclass{53D20}
\secondaryclass{53D45}
\keywords{Gromov width, Moser's method, symplectic embedding, 
complex Grassmannian, moment map}

\maketitle

\section{Introduction}
\label{sec:intro}

Consider the ball of capacity $a$
$$ B(a) = \Big \{ z \in \C^N \ \Big | \ \pi \sum_{i=1}^N |z_i|^2 < a 
          \Big \} , $$
with the standard symplectic form
$\omega_\std = \sum dx_j \wedge dy_j$.
The \textsl{Gromov width} of a $2N$-dimensional symplectic manifold $(M,\omega)$
is the supremum of the set of $a$'s such that $B(a)$ can be symplectically
embedded in $(M,\omega)$.
Computations of Gromov width and, more generally, of symplectic ball
packings, can be found, for example, in \cite{Bi1,Bi2,Bi3,G,MP,T}.

Often in symplectic geometry, equivariant techniques
give constructions whereas holomorphic techniques give obstructions.
We use both.

Our main technical result is a criterion for the existence of symplectic 
embeddings of open subsets of $\C^n$ into a symplectic manifold
with a Hamiltonian torus action.
See Proposition \ref{embed}.   

Paul Biran has asked whether the Gromov width of a compact symplectic 
manifold is at least one if the symplectic form is integral.  
In Proposition \ref{wts 1} we answer his question positively
whenever the manifold admits a Hamiltonian circle action
with a fixed point $p$ such that all the isotropy weights at $p$
are equal to one.  

As a corollary, we obtain embeddings of balls into complex Grassmannians.
More precisely, we prove the following theorem.

\begin{Theorem} \label{thm:gromov}
Let $\Gr(k,n)$ be the Grassmannian of $k$-planes in $\C^n$, 
together with its $\U(n)$-invariant symplectic form $\omega$, 
normalized so that $[\omega]$ generates the integral cohomology 
$H^2(\Gr(k,n);\Z)$.
There exists a symplectic embedding of $B(a)$ into $\Gr(k,n)$
if and only if $a \leq 1$.
\end{Theorem}

Next, we use holomorphic techniques to show that it is 
impossible to embed the ball $B(a)$ into $\Gr(k,n)$ if $a > 1$.
The proof uses two ingredients:
a slight adaptation of the proof of Gromov's non-squeezing theorem,
and the calculation of certain Gromov-Witten invariants 
for $\Gr(k,n)$, which we quote from \cite{ST}
(also see \cite{kim,BDW}).

The Gromov width of the complex Grassmannian was independently computed
by Guangcun Lu in \cite{glu}. Lu obtained the lower bound by an explicit
embedding of a ball.
Our results are more general in that they give lower bounds
for the Gromov width of many more manifolds, 
such as other generalized flag manifolds, and in that they allow one to embed sets 
other than balls.

\subsubsection*{Acknowledgment}
We thank Paul Biran, Kai Cieliebak, Guangcun Lu, Dusa McDuff, 
Ignasi Mundet i Riera, and Martin Pinsonnault for useful discussions.

During work on this project,
Y. Karshon was partially supported by an NSERC Discovery grant,
S. Tolman was partially supported by a Sloan fellowship and by
the National Science Foundation grant DMS \#02-04448,
and both authors were partially supported by the United States -- Israel 
Binational Science Foundation grant number 2000352.
 

\section{Lower bounds for Gromov width}

In this section, we construct symplectic embeddings of open subsets of
$\C^n$ into  symplectic manifolds with  Hamiltonian torus actions.
The key technique is Moser's method.
This section is an extension of our work in
\cite[\S 13]{locun} and is inspired by \cite[\S 1]{De}.

Let  a torus $T \cong (S^1)^{\dim T}$ with Lie algebra $\t$
act effectively on a connected symplectic manifold $(M,\omega)$
by symplectic transformations.
A \textsl{moment map} is a map $\Phi \colon M \to \t^*$ such that
\begin{equation} \label{def moment}
        \iota(\xi_M) \omega = - \d \left< \Phi,\xi \right>
\quad \forall \ \xi \in \t,
\end{equation}
where $\xi_M$ is the corresponding vector field on $M$.

Let $p \in M$ be a fixed point. There exist $\eta_j \in \t^*$,
called the \textsl{isotropy weights} at $p$, such that
the induced linear symplectic $T$-action on the tangent space $T_pM$
is isomorphic to the action on $(\C^n,\omega_\std)$
generated by the moment map
$$ \Phi_{\C^n}(z) = \Phi(p) + \pi \sum |z_j|^2 \eta_j .$$
The isotropy weights are uniquely determined up to permutation.

By the equivariant Darboux theorem \cite{W},
a neighborhood of $p$ in $M$ is equivariantly symplectomorphic
to a neighborhood of $0$ in $\C^n$.
The results of this section allow us to control the size
of this neighborhood.

For the applications in this paper it is enough to symplectically embed
the ball $B(1) \subset \C^n$ into manifolds with $S^1$-actions;
see Proposition \ref{wts 1}.
However, we will take this opportunity to develop the relevant machinery 
in the more general case
where we embed other subsets of $\C^n$, possibly unbounded, 
into manifolds with torus actions.

Let $\calT \subset \ft^*$ be an open convex set which contains $\Phi(M)$.
The quadruple $(M,\omega,\Phi,\calT)$ is a
\textsl{proper Hamiltonian $\mathbf{T}$-manifold} if 
$\Phi$ is proper as a map to $\calT$, that is,
the preimage of every compact subset of $\calT$ is compact.

For any subgroup $K$ of $T$, let
$M^K = \{ m \in M \mid a \cdot m = m \ \forall a \in K \}$
denote its fixed point set.

\begin{Definition} \label{centered-definition}
A proper Hamiltonian $T$-manifold $(M,\omega,\Phi,\calT)$
is \textsl{centered} about a point $\alpha \in \calT$ if
$\alpha$ is contained in the moment map image of every component
of $M^K$, for each $K \subseteq T$.
\end{Definition}

\begin{Example}
A compact symplectic manifold with a non-trivial $T$-action is never
centered, because it has fixed points with different moment map images.
\end{Example}

\begin{Example}
Let a torus $T$ act linearly on $\C^n$ with a proper moment map
$\Phi_{\C^n}$ such that $\Phi_{\C^n}(0) = 0$.  
Let $\calT \subset \t^*$ is be an open convex subset containing the origin. 
Then $\Phi_{\C^n}\inv(\calT)$ is centered about the origin.
\end{Example}

\begin{Example}
Let $M$ be a compact symplectic toric manifold with moment map
$\Phi \colon M \to \t^*$.   Then $\Delta := \image \Phi$
is a convex polytope.  The orbit type strata in $M$ are the
moment map pre-images of the relative interiors of the faces
of $\Delta$.  Hence, for any $\alpha \in \Delta$,
$$ \bigcup\limits_{\substack{F \text{ face of } \Delta \\ \alpha \in F}}
   \Phi\inv(\text{rel-int } F) $$
is the largest subset of $M$ that is centered about $\alpha$.
\end{Example}

\begin{Example}
Let $(M,\omega,\Phi,\calT)$ be a proper Hamiltonian $T$-manifold.
Then every point in $\t^*$ has a neighborhood whose preimage is centered.
This is a consequence of the local normal form theorem and the properness
of the moment map.
\end{Example}

\begin{Remark}
In \cite[Definition 1.4]{locun} we defined a proper Hamiltonian
$T$-manifold $(M,\omega,\Phi,\calT)$ to be centered
about $\alpha \in \calT$ if $\alpha$ is contained in the closure
of the moment map image of every orbit type stratum.
This is equivalent to Definition \ref{centered-definition} above 
because of two facts.
First, the components of the fixed point sets of subgroups of $T$
are precisely the closures of the orbit type strata.
Second, because $\Phi$ is proper, the closure of the image
of any subset is equal to the image of its closure.
\end{Remark}

\begin{Proposition} \label{embed}
Let $(M,\omega,\Phi,\calT)$ be a proper Hamiltonian $T$-manifold.
Assume that $M$ is centered about $\alpha \in \calT$ and that 
$\Phi\inv(\{\alpha\})$ consists of a single fixed point $p$.
Then $M$ is equivariantly symplectomorphic to 
$$\left\{ z \in \C^n \ | \ \alpha + \pi \sum |z_j| \eta_j \in \calT \right\},$$
where $\eta_1,\ldots,\eta_n$ are the isotropy weights at $p$.
\end{Proposition}

\begin{proof}
For simplicity, assume that $\alpha = 0$.  Let
$$ \Phi_{\C^n}(z) = \pi \sum |z_j|^2 \eta_j. $$
There exist a convex neighborhood $V$ of $0$ in $\calT$ and
a $T$-equivariant symplectomorphism from $\Phi\inv(V) \subseteq M$
to $\Phi_{\C^n}\inv(V) \subseteq \C^n$:
$$ F \colon \Phi\inv(V) \to (\Phi_{\C^n})\inv(V) .$$
These exist for the following reasons.
The equivariant Darboux theorem gives an equivariant symplectomorphism 
from a neighborhood of $p$ in $M$ to a neighborhood of $0$ in $\C^n$.
Because $\Phi\inv(\{0\})$ consists of a single point,
the moment map $\Phi_{\C^n} \colon \C^n \to \t^*$ is proper.
(See, for example, Lemma 5.4 in \cite{locun}.)
Hence, the neighborhoods of $p \in M$ and $0 \in \C^n$
contain the preimages of a neighborhood of $0 \in \t^*$.

The \textsl{Euler vector field} on a vector space is the generator
of the flow $x \mapsto e^t x$.
Let $X$ be the Euler vector field on $\t^*$;
then $-X$ generates the flow $ x \mapsto e^{-t} x.$

Since $(M,\omega,\Phi,\calT)$ is a proper Hamiltonian $T$-manifold 
that is centered about the origin $0 \in \t^*$,
by Lemma 13.2 of \cite{locun}\footnote{
For our application in Proposition \ref{wts 1}, 
we only need to consider the case that ${\dim T = 1}$,
and $\Phi \colon M \to \R$ is a submersion away from $p$.
Since $\Phi$ is  homogeneous near $p$ with respect to Darboux coordinates,
we can lift the vector field $X$ to $M$
without referring to Lemma 13.2 of \cite{locun}.
}
there exists a smooth invariant 
vector field $\widetilde{X}$ on $M$ such that 
$\Phi_*(\widetilde{X}) = X$.
The vector field $-\widetilde{X}$ generates a flow
$\widetilde{\psi}_t \colon M \to M,$
which satisfies
$$ 
\Phi \circ \widetilde{\psi}_t = e^{-t}\, \Phi
$$
wherever it is defined.
Hence, since $\Phi$ is proper and $\calT$ is convex
and contains the origin,
$\widetilde{\psi}_t$ is  defined for all $t \geq 0$.
Because $\widetilde{X}$ is $T$-invariant,
$\widetilde{\psi}_t$ is $T$-equivariant.
Define a family of $T$-invariant symplectic forms on $M$
by
$$ \omega_t = e^t \, (\widetilde{\psi}_t)^* \omega .$$
Let $\widehat{X}$ be half the Euler vector field on $\C^n$.
Then $-\widehat{X}$ generates the flow
$\widehat{\psi}_t \colon \C^n \to \C^n$ given by
$$ \widehat{\psi}_t(z) =  e^{-t/2} \, z. $$
Note that  
$\Phi_{\C^n} \circ \widehat{\psi}_t = e^{-t} \, \Phi_{\C^n}$ and
$ (\widehat{\psi}_t)^* \omega_\std = e^{-t} \, \omega_\std$.
So
$$ \widetilde{F}_t := 
(\widehat{\psi}_t)\inv \circ F \circ \widetilde{\psi}_t
\colon (\Phi\inv(e^t V),\omega_t) 
\to (\Phi_{\C^n}\inv (e^t V \cap \calT) , \omega_\std) $$
is an equivariant symplectomorphism.

We will now apply Moser's method.
Let $\lambda$ be a $T$-invariant one form such that $d\lambda = \omega$.
Let
$$ \lambda_t = e^t (\widetilde{\psi}_t)^* \lambda. $$
Then $d\lambda_t = \omega_t$. Also, $\iota(\xi_M) \lambda_t = \Phi^\xi$,
because both sides take the value $0$
at $p$, and they have the same differential.
Hence, $\beta_t := \frac{d}{dt} \lambda_t$ is $T$-invariant
and $\iota(\xi_M) \beta_t = 0$ for all $\xi \in \t$.
Let $Y_t$ be the time-dependent vector field on $M$ which is determined
by $i_{Y_t} \omega_t = -\beta_t$.  Then $Y_t \in \ker d \Phi$ because
$$\left< d\Phi(Y_t), \xi \right>
 = - \omega_t(\xi_M , Y_t) = - i_{\xi_M} \beta_t = 0$$
for all $\xi \in \t$.
Since $\Phi$ is proper, this implies that $Y_t$
integrates to an isotopy, $G_t \colon M \to M$.
Since $\omega_t$ and $\beta_t$ are $T$-invariant,
the vector field $Y_t$ is also $T$-invariant;
consequently, $G_t$ is $T$-equivariant.
Finally,
$$\frac{d}{dt}(G_t^* \omega_t) =
        G_t^* (L_{Y_t} \omega_t) + G_t^* (\frac{d}{dt} \omega_t)
        = G_t^* (-d\beta_t) + G_t^* (d\beta_t) = 0. $$
Hence, $G_t^* \omega_t$ is independent of $t$.
Because $G_0^*\omega_0 = \omega$,
$$ G_t \colon (M,\omega) \to (M,\omega_t) $$
is a symplectomorphism for all $t \geq 0$.

Thus, the composition
$$
 \widetilde{F}_t \circ G_t \colon (\Phi\inv(e^t \, V),\omega) \to 
(\Phi_{\C^n}\inv(e^t \,V \cap \calT),\omega_\std) 
$$
is an equivariant symplectomorphism.
If $\calT$ is bounded, then we can choose
$t$  sufficiently large so that $\calT \subset e^t \, V$;
hence we are done.  

Now suppose that $\calT$ is not bounded. 
We will modify our constructions of the vector field
$\tilde{X}$ and the one-form $\lambda$
so that
$\widetilde{F}_t  = \widetilde{F}_s$ and 
$G_s = G_t$  
on the intersection $\Phi\inv(e^t V) \cap \Phi\inv(e^s V)$.

First, we modify our  construction of $\widetilde{X}$ so that,
in addition to satisfying $\Phi_*(\widetilde{X})=X$,
after possibly shrinking $V$
\begin{equation} \label{property1}
F_*(\widetilde{X}) = \widehat{X} \quad \text{on} \quad \Phi\inv(V).
\end{equation}
To do this, we construct $\widetilde{X}$ on $M \ssminus \{p\}$
as before, and then
patch with the vector field $F^*(\widehat{X})$ on $\Phi\inv(V)$,
using an invariant partition of unity subordinate to the sets 
$M \ssminus \{p\}$ and $\Phi\inv(V)$.

The property \eqref{property1}
implies that $F \circ \widetilde{\psi}_t = \widehat{\psi}_t \circ F$
on $\Phi\inv(V)$, so 
$$
 \widetilde{F}_t = \widetilde{F}_s \quad \text{on} \quad
\Phi\inv(e^t V) \cap \Phi\inv(e^s V) .
$$
Therefore, since the union of the sets $e^tV$ over all $t \geq 0$
is all of $\t^*$,
we can define an equivariant diffeomorphism
$$ \widetilde{F}  \colon M \to \Phi_{\C^n}\inv (\calT) $$
$$ \widetilde{F} = F_t
   \quad \text{on} \quad \Phi\inv(e^tV).\leqno{\hbox{by}}$$
$$\widetilde{F}^*(\omega_\std) = \omega_t \quad
\text{on} \quad \Phi\inv(e^tV).\leqno{\hbox{Clearly,}}$$
Next, we modify our construction of the invariant
one-form $\lambda$ so
that, in addition to satisfying $d \lambda = \omega$,
after possibly shrinking $V$
\begin{equation} \label{property2}
F^* (\lambda_\std) = \lambda \quad \text{on} \quad \Phi\inv(V).
\end{equation}
Consider the one-form 
$\lambda_\std = \half \sum_{i=1}^n (x_i dy_i - y_i dx_i)$
on $\C^n$.
Let $\lambda_V$ be any $T$-invariant one-form on $M$ such that
$\lambda_V = F^* \lambda_\std$ on $\Phi\inv(V)$.
Then $\omega' := \omega - d\lambda_V$ is a closed two-form on $M$
which vanishes on $\Phi\inv(V)$. Since $M$ is diffeomorphic
to $\Phi_{\C^n}\inv(\calT)$, $\calT$ is convex,
and $\Phi_{\C^n}$ is homogeneous, $M$ is contractible;
similarly, $\Phi\inv(V)$ is contractible.
Consequently, the relative cohomology $H^2(M,\Phi\inv(V))$ is zero.
This implies that $\omega'$ has a primitive one-form $\lambda'$
which vanishes on $\Phi\inv(V)$; we may choose it to be $T$-invariant.
Let $\lambda = \lambda' + \lambda_V$. 

Since
$ (\widehat{\psi_t})^* \lambda_\std = e^{-t} \lambda_\std$,
the property \eqref{property2} implies that 
$$ \lambda_t = \widetilde{F}_t^* (\lambda_\std) \quad \text{on} \quad
\Phi\inv(e^t V).$$
Hence, $\lambda_t = \lambda_s$ 
on $\Phi\inv(e^t V) \cap \Phi\inv(e^s V)$.
Therefore $\beta_t = 0$, and hence $Y_t = 0$,
on $\Phi\inv(e^t V)$.  Consequently,
$$
 G_t = G_s \quad \text{on} \quad \Phi\inv(e^t V) \cap \Phi\inv(e^s V) .
$$
Therefore, since the union of the sets $e^tV$ over all $t \geq 0$
is all of $\t^*$,
we can define an equivariant diffeomorphism
$$ G  \colon M \to M $$
$$ G = G_t
   \quad \text{on} \quad \Phi\inv(e^tV).\leqno{\hbox{by}}
$$
$$G^*(\omega_t) = \omega \quad
\text{on} \quad \Phi\inv(e^tV).\leqno{\hbox{Clearly,}}$$
Hence, 
$\widetilde{F} \circ G: 
(M,\omega) \to (\Phi_{\C^n}\inv(\calT), \omega_\std)$ is  
an equivariant symplectomorphism, as required. 
\end{proof}

Proposition \ref{embed} answer Biran's question affirmatively 
in a special case:

\begin{Proposition} \label{wts 1}
Let $N$ be a compact 
manifold with an integral symplectic form $\omega$.
Suppose that it admits a Hamiltonian circle action
with a fixed point $p$ such that all the isotropy weights at $p$
are equal to one.  Then there exists a symplectic embedding
of the ball $(B(1),\omega_\std)$ in $(N,\omega)$.
\end{Proposition}

\begin{proof}
We assume, without loss of generality, that $N$ is connected.
Let $\xi_N$ denote the vector field that generates 
the circle action.
Our convention is that the circle group is $S^1 = \R / \Z$, 
so that $\xi_N$ generates a flow of period one.  Let
$$ \Phi \colon N \to \R $$
be the moment map, so that $\iota(\xi_N) \omega = - d \Phi$.
For simplicity, assume that $\Phi(p) = 0$

Because the isotropy weights are positive,
$p$ is an isolated local minimum for the moment map.
Since the moment map fibers are connected \cite{At, GS},
$\Phi\inv(\{0\}) = \{ p \}$.

By Stokes's theorem,
for any fixed point $q$, the difference $\Phi(q) - \Phi(p)$
is equal to the integral of $\omega$ over the cycle obtained
from a curve connecting $p$ to $q$ by ``sweeping" the curve
by the circle action.
Because $[\omega]$ is integral, this implies that
$\Phi(q) - \Phi(p)$ is an integer.
So $p$ is the only fixed point that is mapped to $[0,1)$.

Let $M = \Phi\inv([0,1))$.
Consider a subgroup $K \subset S^1$ and let $Y \subset M^K$
be a connected component of its fixed point set.
Since $Y$ is closed in $M$ and $\Phi\colon M \to [0,1)$ 
is proper,
the image of $Y$ is a closed subset of $[0,1)$,
so it has a minimum.  Any point in $Y$ which is mapped
to this minimum must be a fixed point.  Hence $p \in Y$,
and so $0 \in \Phi(Y)$.  This shows that $M$ is centered
about $0$.

Proposition \ref{wts 1} then follows from Proposition \ref{embed}.
\end{proof}

\section{Lower bounds for Grassmannians}
\label{sec:lower for orbits}

In this section we construct
an embedding of the ball $B(1)$ into the complex Grassmannian $\Gr(k,n)$,
thus showing that the Gromov width of the Grassmannian is at least one.

\begin{Proposition} \label{lower for Gkn}
Let $\Gr(k,n)$ be the Grassmannian of $k$-planes in $\C^n$,
together with its $\U(n)$-invariant symplectic form $\omega$,
normalized so that $[\omega]$ generates the integral cohomology
$H^2(\Gr(k,n);\Z)$.  There exists a symplectic embedding of $B(1)$
into $\Gr(k,n)$.
\end{Proposition}

\begin{proof}
Let the circle group $S^1$ act on $\C^n$ by
$$ a \cdot (z_1,\ldots,z_n)
  = (z_1, \ldots, z_k, a z_{k+1}, \ldots, a z_n).$$
Take the induced action on $\Gr(k,n)$.
Then $p = \C^k \times \{ 0 \}$ is a fixed point for this action.
Since $T_p \Gr(k,n) \cong \text{Hom}(\C^k,\C^{n-k})$,
the isotropy action is complex multiplication by $a \in S^1$.
Proposition \ref{lower for Gkn} now follows
from Proposition~\ref{wts 1}.
\end{proof}

\section{Upper bounds for Gromov width}

In this section, we give a short review of Gromov-Witten invariants
and how they can be used to give upper bounds to Gromov widths.
This material appears in detail in many places; our treatment is adapted 
from \cite{MS}.

Let $(M,\omega)$ be a compact symplectic manifold.
A homology class $B \in H_2(M)$ is \textsl{spherical}
if it is in the image of the Hurewicz homomorphism
$\pi_2(M) \to H_2(M)$.
A homology class $B \in H_2(M) $ is \textsl{indecomposable}
if it does not decompose as a sum $B = B_1 + \cdots + B_k$ 
of spherical classes such that $\omega(B_i) > 0$.

The form $\omega$ \textsl{tames}
an  almost complex structure $J$ on $M$ if 
$$\omega(v,Jv) >0 $$
for all non-zero $v \in TM$. Given $\omega$,
there are many almost complex structures $J$ on $M$
that are tamed by $\omega$; however, the Chern
classes $c_i(TM)$ of the complex vector bundles $(TM,J)$ 
are independent of this choice.
We say that $(M, \omega)$ is \textsl{monotone} if 
there exists a positive constant $\lambda > 0$ so that
$\omega(B) = \lambda c_1(TM)(B)$
for all spherical classes $B \in H_2(M;\Z)$.

Assume that $(M,\omega)$ is monotone of dimension $2n$.
Fix an indecomposable spherical class $A \in H_2(M;\Z)$.
Let $N_i$, for $i=1,\ldots,s$, be submanifolds such that 
$\sum_{i=1}^s \codim N_i = 2n + 2c_1(TM)(A) + 2s - 6$.
Let $B_i$ be the homology class represented by $N_i$.
The Gromov invariant\footnote{This invariant can be defined
under more general assumptions.}
$$ \Phi_A(B_1,\ldots,B_s) \in \Z,$$
which is defined, for example, in \cite{MS}, has the following property.
Let $\calJ$ denote the space of almost complex structures
that are tamed by $\omega$.
Let $\calJ_\reg(A) \subseteq \calJ$  denote the set of 
regular almost complex structures for the class $A$ (see \cite{MS}).
Given $J \in \calJ_\reg(A)$, 
for generic deformations $N_i'$ of $N_i$
and generic points $t_i \in \CP^1$,
the number of $J$-holomorphic maps $\CP^1 \to M$ in the class $A$
which send each $t_i$ into $N_i'$, counted with appropriate signs,
is equal to $\Phi_A(B_1,\ldots,B_s)$.  
In particular, if $[p]$ is the homology class of a point
and $X,Y$ are submanifolds such that $\Phi_A([p],[X],[Y]) \neq 0$, 
then, for $J \in \calJ_\reg(A)$,
for every point $p \in M$ and every neighborhood $U$ of $p$ 
there exists a $J$-holomorphic sphere in the class $A$ 
that passes through $U$.

\begin{Proposition} \label{thm:width}
Let $(M,\omega)$ be a monotone symplectic manifold.
Let $A \in H_2(M;\Z)$ be an indecomposable spherical class.
Let $$ \lambda = \int_A \omega.$$
Let $[p]$ denote the homology class of a point.
Suppose that there exist submanifolds $X$ and $Y$ of $M$ so that 
$\dim X + \dim Y = 4n - 2c_1(TM)(A)$
and so that
$$\Phi_A([p],[X],[Y]) \neq 0.$$
If $a > \lambda$, there does not exist 
a symplectic embedding $B(a) \to M$.
\end{Proposition}

\begin{proof}
Suppose that there exists a symplectic embedding 
$$ B(a) \stackrel{\rho}{\hookrightarrow} M $$
where $a > \lambda$.
The standard complex structure $J_\std$ on $B(a) \subset \C^n$
transports through $\rho$ to a complex structure on $\rho(B(a))$.
By a well known technique,
after passing to a smaller $a$ such that $a > \lambda$ there exists
a tamed almost complex structure $J$ on $M$ such that
$\rho$ intertwines $J_\std$ with $J$.
(Since Riemannian metrics can be patched together by partitions of unity,
this follows from the existence of an $\text{Sp}(\R^{2n})$-equivariant projection
from the space of inner products on $\R^{2n}$ to the subspace 
of those inner products that are compatible with the symplectic form.)

Because $\calJ_\reg(A)$ is of second category in $\calJ$,
and by the characterization of Gromov-Witten invariants discussed above,
there exists a sequence $J_i$
of almost complex structures converging to $J$,
and points $p_i$ which converge to $p := \rho(0)$,
and, for each $i$, a $J_i$-holomorphic curves $\CP^1 \to M$ 
in the class $A$ which passes through $p_i$.

Since the class $A$ is indecomposable, by Gromov's compactness theorem
\cite{G,MS} there exists a subsequence of the $J_i$'s which converges 
weakly to a $J$-holomorphic curve $C$ in the class $A$.  
In particular, $\int_C \omega = \lambda$.
Since $p_i$ converge to $p$, this limit curve contains $p$.

The pre-image of $C$ under $\rho$ is a holomorphic curve in the ball $B(a)$
which passes through the center of the ball
and which is closed in $B(a)$.
But the smallest area of a such a curve is that of a disk through 
the center, which is $a$.  (See \cite[page 99]{AL}.)
This contradicts the assumption that $a > \lambda$.
\end{proof}

\section{Upper bounds for Grassmannians}

It remains to prove the following proposition.

\begin{Proposition} \label{upper for Gkn}
Let $\Gr(k,n)$ be the Grassmannian of $k$-planes in $\C^n$,
together with its natural $\U(n)$ invariant symplectic form,
normalized so that $[\omega]$ generates the integral cohomology
$H^2(\Gr(k,n);\Z)$.
If $a > 1$, there does not exist a symplectic embedding 
of $B(a)$ into $\Gr(k,n)$.
\end{Proposition}

\begin{proof}
The real dimension of $M = \Gr(k,n)$ is $2k(n-k)$.
Let $A$ be the generator of $H_2(\Gr(k,n),\Z) \cong \Z$.
Clearly, $A$ is irreducible.

The standard complex structure is tamed by $\omega$, and
$c_1(TM)(A) = n$.  It is easy to check that $\Gr(k,n)$ is monotone.

Fix a hyperplane $W \subset \C^n$ and let $X \subset \Gr(k,n)$ be
the set of $k$-planes that are contained in $W$.
Fix a vector $y \in \C^n$ and let $Y \subset \Gr(k,n)$
be the set of $k$-planes that contain $y$.
Let $p \in \Gr(k,n)$ be any point.

These submanifolds represent homology classes
$[p]$, $[X]$, and $[Y]$ in $H_0(M)$, $H_{2k(n-1)}(M)$, and 
$H_{2n(k-1)}(M)$, respectively.

In \cite{ST}, it was shown that
$\Phi_A([p],[X],[Y]) = 1$.
The Proposition then follows from Proposition \ref{thm:width}.
\end{proof}

Theorem \ref{thm:gromov} follows immediately from Propositions
\ref{lower for Gkn} and \ref{upper for Gkn}.

\Addresses\recd

\end{document}